\documentclass[12pt]{article}
\usepackage{a4}
\usepackage[fleqn]{amsmath}
\usepackage{amsthm}
\usepackage{amssymb}
\usepackage{xypic}

\newtheorem{proposition}{Proposition}
\newtheorem{corollary}[proposition]{Corollary}
\newtheorem{theorem}[proposition]{Theorem}
\newtheorem{lemma}[proposition]{Lemma}

\theoremstyle{definition}

\newtheorem{example}{Example}
\newtheorem{remark}{Remark}

\newcommand{\B}[1]{\mathbf{#1}}

\newcommand{\g}{\mathfrak g}
\newcommand{\h}{\mathfrak h}
\newcommand{\n}{\mathfrak{n}}

\renewcommand{\sl}{\mathfrak{sl}}

\DeclareMathOperator{\id}{id}

\DeclareMathOperator{\ad}{ad} \DeclareMathOperator{\Tr}{Tr}
\DeclareMathOperator{\Ker}{Ker} \DeclareMathOperator{\Image}{Im}

\DeclareMathOperator{\End}{End} \DeclareMathOperator{\spanv}{span}

 \DeclareMathOperator{\Fun}{Fun}
\DeclareMathOperator{\Hom}{Hom} \DeclareMathOperator{\fin}{fin}
\DeclareMathOperator{\rrr}{r} 
 \DeclareMathOperator{\red}{red}
\DeclareMathOperator{\Mat}{Mat}

\def\ftext#1{{\let\thefootnote\relax
\footnotetext{\vskip-.9\baselineskip\noindent #1}}}

\begin{document}
\title{Irreducible highest-weight modules and equivariant quantization}
\author{E.~Karolinsky$^{\,\diamond}$, A.~Stolin,
and V.~Tarasov$^{\,*}$}
\date{}

\maketitle

\ftext{$^{\diamond}$Supported in part by the Royal Swedish Academy of Sciences
\\
$^*$Supported in part by RFFI grant 02--01--00085a and
CRDF grant RM1--2334--MO--02}

\begin{sloppy}

\vskip-.6\baselineskip
\hrule height0pt
\thispagestyle{empty}

\section{Introduction}

The notion of deformation quantization, motivated by ideas coming
from both physics and mathematics, was introduced in classical
papers \cite{BFFLS, FLS1, FLS2}. Roughly speaking, a deformation
quantization of a Poisson manifold $(P,\{ \,, \})$ is a formal
associative product on $(\Fun P)[[\hbar]]$ given by $f_1\star
f_2=f_1f_2+\hbar c(f_1,f_2)+O(\hbar^2)$ for any $f_1, f_2\in \Fun
P$, where the skew-symmetric part of $c$ is equal to $\{ \,, \}$,
and the coefficients of the series for $f_1\star f_2$ should be
given by bi-differential operators.

The fact that any Poisson manifold can be quantized in this sense
was proved by Kontzevich in \cite{Kontz1}. However, finding exact
formulas for specific cases of Poisson brackets is an interesting
separate problem. There are several well-known examples of such
explicit formulas. One of the first was the Moyal product
quantizing the standard symplectic structure on $\mathbb R^{2n}$.
Another one is the standard quantization of the
Kirilov-Kostant-Souriau bracket on the dual space $\g^*$ to a Lie
algebra $\g$ (see \cite{Gutt}). Relations between this
quantization and the Yang-Baxter equation was shown by Gekhtman
and Stolin in \cite{GS}.

Despite the formula for the standard quantization of the
Kirilov-Kostant-Souriau bracket is known already for a long time,
the problem of finding explicit formulas for equivariant
quantization of its symplectic leaves, i.e., coadjoint orbits on
$\g^*$, was open. Recently this problem was solved in important
cases in \cite{AL, DM, EEM, EKKMASVT} using the relationship with
the dynamical Yang-Baxter equation and the Shapovalov form on
Verma modules.

This paper is a continuation of \cite{EKKMASVT}. One of the main
results obtained in \cite{EKKMASVT} is the connection between
quantum dynamical twists and equivariant quantization. More
precisely, let $\g=\n_-\oplus\h\oplus\n_+$ be a finite-dimensional
complex semisimple Lie algebra with a fixed triangular
decomposition, $F={\mathbb C}[G]$ the algebra of all matrix
elements of all finite dimensional representations of $\g$,
$M(\lambda)$ the Verma module with highest weight
$\lambda\in\h^*$, $J(\lambda)$ the universal fusion element that
corresponds to $\lambda$. Assume that $\lambda$ is generic, i.e.,
$M(\lambda)$ is irreducible. We have a natural map
$\Hom_\g(M(\lambda),M(\lambda)\otimes F)\to F[0]$ which is an
isomorphism of vector spaces. The translation onto $F[0]$ of the
natural product on $\Hom_\g(M(\lambda),M(\lambda)\otimes F)$ is
given by the formula $f_1 \star_\lambda
f_2=\mu\left(\overrightarrow{J(\lambda)}(f_1\otimes f_2)\right)$,
where $\mu$ is the initial product on $F$ restricted onto $F[0]$.
We may treat the obtained algebra as an equivariant quantization
of the coadjoint orbit $\mathcal O_\lambda\subset\g^*$ equipped
with the Kirillov-Kostant-Souriau bracket.

The main goal of this paper is to present some generalizations of
the above mentioned results to the case of non-generic $\lambda$.
In fact, we give explicit formulas for star-products on certain
subspaces of $F[0]$, which are in general not closed under the
original multiplication on $F$.

Consider the irreducible $\g$-module $V(\lambda)$ with highest
weight $\lambda\in\h^*$. We have
$V(\lambda)=M(\lambda)/K_\lambda\mathbf{1}_\lambda$, where
$K_\lambda\subset U\n_-$, and $\mathbf{1}_\lambda$ is the
generator of $M(\lambda)$. Consider also the opposite Verma module
$\widetilde{M}(-\lambda)$ with the lowest weight $-\lambda\in\h^*$
and the lowest weight vector $\widetilde{{\mathbf 1}}_{-\lambda}$.
Note that the maximal $\g$-submodule in $\widetilde{M}(-\lambda)$
is of the form $\widetilde{K}_{\lambda}\cdot\widetilde{{\mathbf
1}}_{-\lambda}$, where $\widetilde{K}_{\lambda}\subset U\n_+$. We
get a vector space isomorphism
$\Hom_\g(V(\lambda),V(\lambda)\otimes F)\simeq
F[0]^{K_\lambda+\widetilde{K}_\lambda}$, which allows one to
consider the product on $F[0]^{K_\lambda+\widetilde{K}_\lambda}$
induced by the natural multiplication in
$\Hom_\g(V(\lambda),V(\lambda)\otimes F)$. One can express this
product on $F[0]^{K_\lambda+\widetilde{K}_\lambda}$ as
$f_1\star_\lambda
f_2=\mu\left(\overrightarrow{J_{\red}(\lambda)}(f_1\otimes
f_2)\right)$, where the ``reduced fusion element''
$J_{\red}(\lambda)$ can be computed in terms of the Shapovalov
form on $V(\lambda)$.

For a special case when $\lambda\in\h^*$ satisfies
$\langle\lambda,\alpha^\vee\rangle=0$ for any $\alpha$ in some
simple root subset $\Delta$ and generic otherwise we show that
$F[0]^{K_\lambda+\widetilde{K}_\lambda}\simeq F[0]^{K_\lambda}$.
In this case $F[0]^{K_\lambda}$ is closed under the original
multiplication on $F$. In fact, this is the algebra of regular
functions on the coadjoint orbit $\mathcal O_\lambda$. Hence the
algebra $(F[0]^{K_\lambda},\star_\lambda)$ can be viewed as an
equivariant quantization of $\mathcal O_\lambda$.

Finally, we investigate limiting properties of the universal
fusion element $J(\lambda)$. In particular we show that for some
values of $\lambda_0\in\h^*$ we can guarantee that
$f_1\star_\lambda f_2\rightarrow f_1\star_{\lambda_0}f_2$ as
$\lambda\rightarrow\lambda_0$. We also show that for any
$\lambda_0$ having a ``good limiting property'' of this type the
action map $U\g\to(\End V(\lambda_0))^{\rrr}_{\fin}$ is surjective
(here $(\End V(\lambda_0))^{\rrr}_{\fin}$ stands for the locally
finite part of $\End V(\lambda)$ with respect to the adjoint
action of $U\g$). Note that this surjectivity question is known as
the classical problem of Kostant (see \cite{Ja, Jo}). The complete
answer to this question is still unknown. However, there are
examples of $\lambda_0$ such that $U\g\to(\End
V(\lambda_0))^{\rrr}_{\fin}$ is not surjective (see \cite{Jo}).
There is also a known class of simple highest weight modules for
which this map is surjective. We comment on the Kostant problem in
other parts of the paper as well.

We also notice that most of the results of this paper have
analogues for quantized universal enveloping algebras. We will
discuss these questions in details elsewhere.

This paper is organized as follows. Section \ref{Sect-Hopf}
contains some general Hopf-algebraic constructions that will be
useful in the sequel. In Section \ref{Sect-main} we provide a
construction of a star-product on
$F[0]^{K_\lambda+\widetilde{K}_\lambda}$ by means of the
Shapovalov form on $V(\lambda)$. Subsection
\ref{Subsect-main-general} is devoted to the general construction,
and in Subsection \ref{Subsect-main-symm} we discuss applications
to symmetric spaces and coadjoint orbits. Finally, Section
\ref{Sect-limit} is devoted to study of limiting properties of
fusion elements and corresponding star-products.

Throughout this paper all Lie algebras are assumed to be
finite-dimensional, and the ground field is $\mathbb C$.

\subsection*{Acknowledgements} The authors are grateful to Maria
Gorelik and Volodymyr Mazorchuk for useful discussions on the
topic of the paper. Part of this research was done during E.K.'s
visit to University of G\"oteborg; we thank our colleagues there
for their hospitality.

\section{Hopf algebra preliminaries}\label{Sect-Hopf}

Let $A$ be a Hopf algebra. As usual, we will denote by $\Delta$
(resp.\ $\varepsilon$, $S$) the comultiplication (resp.\ counit,
antipode) in $A$. We will systematically use the Sweedler notation
for comultiplication, i.e., $\Delta(x)=\sum_{(x)}x_{(1)}\otimes
x_{(2)}$,
$(\Delta\otimes\id)\Delta(x)=(\id\otimes\Delta)\Delta(x)=\sum_{(x)}x_{(1)}\otimes
x_{(2)}\otimes x_{(3)}$, etc.

Assume $M$ is a (left) $A$-module. An element $m\in M$ is called
{\it locally finite} if $\dim Am<\infty$. Denote by $M_{\fin}$ the
subset of all locally finite elements in $M$. Clearly, $M_{\fin}$
is a submodule in $M$. Similarly, we can consider locally finite
elements in a right $A$-module $N$. For convenience, we will use
the notation $N^{\rrr}_{\fin}$ for the submodule of all locally
finite elements in this case.

Recall that the left (resp.\ right) adjoint action of $A$ on
itself is defined by the formula
$\ad_xa=\sum_{(x)}x_{(1)}aS(x_{(2)})$ (resp.\
$\ad^{\rrr}_xa=\sum_{(x)}S(x_{(1)})ax_{(2)}$). We denote by
$A_{\fin}$ (resp.\ $A^{\rrr}_{\fin}$) the corresponding submodules
of locally finite elements. Since
$\ad_x(ab)=\sum_{(x)}\ad_{x_{(1)}}(a)\ad_{x_{(2)}}(b)$, we see
that $A_{\fin}$ is a (unital) subalgebra in $A$; the same holds
for $A^{\rrr}_{\fin}$. If the antipode $S$ is invertible, then $S$
defines an isomorphism between $A_{\fin}$ and $A^{\rrr}_{\fin}$.
We will assume that $S$ is invertible.

Fix a Hopf subalgebra $F$ of the Hopf algebra $A^\star$ dual to
$A$.
In the sequel we will use the left and right regular actions of
$A$ on $F$ defined respectively by the formulas
$(\overrightarrow{a}f)(x)=f(xa)$ and
$(f\overleftarrow{a})(x)=f(ax)$.

Now let $M$ be a (left) $A$-module. Equip $F$ with the left
regular $A$-action and consider the space $\Hom_A(M,M\otimes F)$.
For any $\varphi,\psi\in\Hom_A(M,M\otimes F)$ define
\begin{equation}\label{eqn-basic-star-prod}
\varphi \ast \psi =(\id\otimes \mu)\circ (\varphi \otimes
\id)\circ \psi,
\end{equation}
where $\mu$ is the multiplication in $F$.
It is straightforward to verify that
$\varphi*\psi\in\Hom_A(M,M\otimes F)$, and this definition equips
$\Hom_A(M,M\otimes F)$ with a unital associative algebra
structure.

Consider the map $\Phi:\Hom_A(M,M\otimes F)\to\End M$,
$\varphi\mapsto u_\varphi$, defined by
$u_\varphi(m)=(\id\otimes\varepsilon)(\varphi(m))$; here
$\varepsilon(f)=f(1)$ is the counit in $F$. In other words, if
$\varphi(m)=\sum_im_i\otimes f_i$, then
$u_\varphi(m)=\sum_if_i(1)m_i$. Using the fact that $\varepsilon$
is an algebra homomorphism it is easy to show that $\Phi$ is an
algebra homomorphism as well.

\begin{lemma}
The map $\Phi$ embeds $\Hom_A(M,M\otimes F)$ into $\End M$.
\end{lemma}

\begin{proof}
If $\varphi\in\Hom_A(M,M\otimes F)$, $\varphi(m)=\sum_im_i\otimes
f_i$, then
\[
\varphi(am)=a\varphi(m)=\sum_i\sum_{(a)}a_{(1)}m_i\otimes
\overrightarrow{a_{(2)}}f_i,
\]
and
\[
u_\varphi(am)=\sum_i\sum_{(a)}(\overrightarrow{a_{(2)}}f_i)(1)a_{(1)}m_i=
\sum_{(a)}a_{(1)}\left(\sum_if(a_{(2)})m_i\right).
\]
Assume now that $u_\varphi=0$, i.e.,
$\sum_{(a)}a_{(1)}\left(\sum_if(a_{(2)})m_i\right)=0$ for any
$a\in A$ and $m\in M$. Then, in particular,
\begin{gather*}
0=\sum_{(a)}S(a_{(1)})a_{(2)}\left(\sum_if(a_{(3)})m_i\right)=\\
\sum_{(a)}\varepsilon(a_{(1)})\left(\sum_if(a_{(2)})m_i\right)=\sum_if_i(a)m_i
\end{gather*}
for any $a\in A$ and $m\in M$. Obviously, this means that
$\varphi=0$.
\end{proof}

From now on we assume that $F$ contains all matrix elements of the
(left) adjoint action of $A$ on $A_{\fin}$. Since $F$ is closed
under the antipode $(Sf)(x)=f(S(x))$, we see that this assumption
is equivalent to the fact that $F$ contains all matrix elements of
the right adjoint action of $A$ on $A^{\rrr}_{\fin}$.

Let $a\in A^{\rrr}_{\fin}$, i.e., for any $x\in A$ we have
$\ad^{\rrr}_xa=\sum_if_i(x)a_i$, where $f_i\in A^\star$, $a_i\in
A$. In fact, we see that $f_i\in F$ by the assumption above.
Define a linear map $\varphi_a:M\to M\otimes F$ by the formula
$\varphi_a(m)=\sum_ia_im\otimes f_i$. Clearly, $\varphi_a$ is well
defined.

\begin{lemma}
For any $a\in A^{\rrr}_{\fin}$ we have
$\varphi_a\in\Hom_A(M,M\otimes F)$.
\end{lemma}

\begin{proof}
Let $b\in A$. Notice that
\[
\sum_{(b)}b_{(1)}\ad^{\rrr}_{b_{(2)}}y=\sum_{(b)}b_{(1)}S(b_{(2)})yb_{(3)}=
y\sum_{(b)}\varepsilon(b_{(1)})b_{(2)}=yb
\]
for any $y\in A$. Therefore for any $x\in A$ we have
\begin{gather*}
\sum_if_i(x)a_ib=(\ad^{\rrr}_xa)b=\sum_{(b)}b_{(1)}\ad^{\rrr}_{b_{(2)}}\ad^{\rrr}_xa=
\sum_{(b)}b_{(1)}\ad^{\rrr}_{xb_{(2)}}a=\\
\sum_{(b)}b_{(1)}\left(\sum_if_i(xb_{(2)})a_i\right)=
\sum_{(b)}\sum_i(\overrightarrow{b_{(2)}}f_i)(x)b_{(1)}a_i,
\end{gather*}
and
\[
\varphi_a(bm)=\sum_ia_ibm\otimes f_i=
\sum_{(b)}\sum_ib_{(1)}a_im\otimes\overrightarrow{b_{(2)}}f_i=b\varphi_a(m).
\]
\end{proof}

Denote by $\Psi:A^{\rrr}_{\fin}\to\Hom_A(M,M\otimes F)$ the linear
map constructed above (i.e., $\Psi:a\mapsto\varphi_a$).

\begin{lemma}
The map $\Psi$ is an algebra homomorphism.
\end{lemma}

\begin{proof}
Let $a,b\in A^{\rrr}_{\fin}$, $x\in A$,
$\ad^{\rrr}_xa=\sum_if_i(x)a_i$, $\ad^{\rrr}_xb=\sum_jg_j(x)b_j$.
Then
\begin{gather*}
\ad^{\rrr}_x(ab)=\sum_{(x)}\ad_{x_{(1)}}(a)\ad_{x_{(2)}}(b)=\\
\sum_{i,j}\sum_{(x)}f_i(x_{(1)})g_j(x_{(2)})a_ib_j=\sum_{i,j}(f_ig_j)(x)a_ib_j.
\end{gather*}
Thus
\[
\varphi_{ab}(m)=\sum_{i,j}a_ib_jm\otimes f_ig_j=
(\varphi_a*\varphi_b)(m)
\]
for any $m\in M$.
\end{proof}

\begin{remark}
It follows directly from the definitions that the composition
$\Phi\Psi$ equals the restriction to $A^{\rrr}_{\fin}$ of the
canonical homomorphism $A\to\End M$, $a\mapsto a_M$.
\end{remark}

Now consider $A^{\rrr}_{\fin}$, $\Hom_A(M,M\otimes F)$ and $\End
M$ as right $A$-modules: $A^{\rrr}_{\fin}$ via right adjoint
action, $\Hom_A(M,M\otimes F)$ via right regular action on $F$
(i.e., $(\varphi\cdot a)(m)=(\id\otimes
\overleftarrow{a})(\varphi(m))$), and $\End M$ in a standard way
(i.e., $u\cdot a=\sum_{(a)}S(a_{(1)})_Mu{a_{(2)}}_M$). Note that
$A^{\rrr}_{\fin}$, $\Hom_A(M,M\otimes F)$ and $\End M$ equipped
with these structures are indeed right $A$-module algebras, i.e.,
the multiplication map is a module morphism, and the unit is
invariant.

\begin{lemma}
The maps $\Phi$ and $\Psi$ are morphisms of right $A$-modules.
\end{lemma}

\begin{proof}
Straightforward.
\end{proof}

\begin{corollary}\label{cor-Phi-Psi-q}
We have the following morphisms of right $A$-module algebras:
\[
A^{\rrr}_{\fin}\stackrel{\Psi}{\longrightarrow}\Hom_A(M,M\otimes
F)^{\rrr}_{\fin}\stackrel{\Phi}{\longrightarrow}(\End
M)^{\rrr}_{\fin},
\]
and $\Phi\Psi$ is the restriction of the canonical morphism $A\to
\End M$. \qed
\end{corollary}

Now let us assume that $F$ contains all matrix elements of the
canonical right $A$-action on $(\End M)^{\rrr}_{\fin}$ (in
particular, it is enough to require that $F$ contains all matrix
elements of all finite dimensional representations of $A$).

\begin{proposition}\label{prop-Hom-End-isom}
The map
\[
\Phi: \Hom_A (M, M\otimes F)^{\rrr}_{\fin}\longrightarrow(\End
M)^{\rrr}_{\fin}
\]
is an isomorphism of right $A$-module algebras.
\end{proposition}

\begin{proof}
We already know that $\Phi$ is an embedding and homomorphism of
right $A$-module algebras. Now let $u \in (\End M)^{\rrr}_{\fin}$.
Then $u\cdot x=\sum^N_{i=1}f_i(x)u_i$, where $f_i \in F$ and $u_i
\in (\End M)^{\rrr}_{\fin}$. We define $\Xi: (\End
M)^{\rrr}_{\fin}\to\Hom_A(M, M\otimes F)$ by the formula
$\Xi(u)(m)=\sum^N_{i=1}u_i(m)\otimes f_i$. It is straightforward
to verify that $\Xi$ is a morphism of right $A$-module algebras.
Therefore the image of $\Xi$ lies in $\Hom_A(M, M\otimes
F)^{\rrr}_{\fin}$. Since $u=\sum^N_{i=1} f_i(1)u_i$, we conclude
that $\Phi\Xi=\id$. Thus $\Phi$ is surjective and it follows that
$\Phi$ is an isomorphism.
\end{proof}

Suppose that the canonical map $A^{\rrr}_{\fin}\to(\End
M)^{\rrr}_{\fin}$ is an epimorphism.

\begin{proposition}\label{prop-Hopf-submod}
Let $N$ be a submodule of $M$. Then $u(N)\subset N$ for any
$u\in(\End M)^{\rrr}_{\fin}$ and $\varphi(N)\subset N\otimes F$
for any $\varphi\in\Hom_A(M, M\otimes F)^{\rrr}_{\fin}$.
\end{proposition}

\begin{proof}
In this case there exists $a\in A^{\rrr}_{\fin}$ such that
$u(m)=am$ for any $m\in M$. Hence $u(n)=an\in N$ for any $n\in N$.
The second statement follows now from Proposition
\ref{prop-Hom-End-isom}.
\end{proof}

\section{Irreducible highest weight modules and equi\-variant quantization for non-generic
$\lambda$}\label{Sect-main}

\subsection{General construction}\label{Subsect-main-general}

Let $\g$ be a finite-dimensional complex semisimple Lie algebra,
$\h$ its Cartan subalgebra. Fix a triangular decomposition
\begin{equation}\label{eqn-triangular-decomp}
\g=\n_+\oplus\h\oplus\n_-.
\end{equation}
Let $\B R$ be the root system of
$\g$ with respect to $\h$, $\Pi$ the set of simple roots that
corresponds to \eqref{eqn-triangular-decomp}, and $\B R_+$ the
corresponding set of positive roots. We denote by $\rho$ the sum
of fundamental weights. For any $\alpha\in\B R_+$ fix non-zero
elements $X_\alpha\in\g_\alpha$ and $Y_\alpha\in\g_{-\alpha}$.

For any
$\lambda\in\h^*$ let
$M(\lambda)$ be the Verma module with the highest weight
$\lambda$ and the highest weight vector ${\mathbf
1}_\lambda$.

Let $U\g$ be the universal enveloping algebra of $\g$ equipped
with the standard Hopf algebra structure. Clearly,
$(U\g)^{\rrr}_{\fin}=U\g$ and it is well known that the canonical
map $U\g\rightarrow (\End M(\lambda))^{\rrr}_{\fin}$ is epimorphic
for any $\lambda \in\h^*$.

Let $F={\mathbb C}[G]$, that is, $F$ consists of all matrix
elements of all finite dimensional representations of $U\g$. Then
$\Hom_\g (M(\lambda), M(\lambda) \otimes F)^{\rrr}_{\fin}= \Hom_\g
(M(\lambda), M(\lambda) \otimes F)$.

Let $K(\lambda)$ be the maximal $\g$-submodule of $M(\lambda)$ and
$V(\lambda)=M(\lambda)/K(\lambda)$ be the irreducible $\g$-module
with highest weight $\lambda \in \h^*$. Applying Proposition
\ref{prop-Hopf-submod} we get the canonical maps $(\End
M(\lambda))^{\rrr}_{\fin}\to (\End V(\lambda))^{\rrr}_{\fin}$ and
$\Hom_\g (M(\lambda), M(\lambda)\otimes F) \to \Hom_\g
(V(\lambda), V(\lambda)\otimes F)$.

We have the following

\begin{proposition}\label{prop-Hom-End-equiv}
Let $\Phi_M$ be the map from Proposition \ref{prop-Hom-End-isom}.
Then the diagram
\[
 \xymatrix{\Hom_\g (M(\lambda), M(\lambda)\otimes F)
\ar[d]^{\Phi_{M(\lambda)}} \ar[r] & \Hom_\g
(V(\lambda), V(\lambda)\otimes F) \ar[d]^{\Phi_{V(\lambda)}} \\
(\End M(\lambda))^{\rrr}_{\fin} \ar[r] & (\End
V(\lambda))^{\rrr}_{\fin}}
\]
is commutative.\qed
\end{proposition}

Denote by $\overline{\mathbf1}_\lambda$ the image of
$\mathbf1_\lambda$ in $V(\lambda)$. For any
$\varphi\in\Hom_\g(V(\lambda),V(\lambda)\otimes F)$ the formula
$\varphi(\overline{\mathbf 1}_\lambda) =\overline{\mathbf
1}_\lambda \otimes f_\varphi + \sum_{\mu<\lambda} v_\mu \otimes
f_\mu$ defines a map $\Theta:\Hom_\g(V(\lambda),V(\lambda)\otimes
F)\to F[0]$, $\varphi\mapsto f_\varphi$.

\begin{theorem}\label{thm-theta-inj}
$\Theta$ is an embedding.
\end{theorem}


We want also to describe the image of $\Theta$. We will need some
extra notation.

Denote by $x\mapsto(x)_0$ the projection $U\g\to U\h$ along
$\n_-\cdot U\g+U\g\cdot\n_+$. For any $\lambda\in\h^*$ consider a
pairing $\pi_\lambda:U\n_+\otimes U\n_-\to\mathbb C$ defined by
$\pi_\lambda(x\otimes y)=(\overline{x}y)_0(\lambda)$ (here
$S:x\mapsto\overline{x}$ is the antipode in $U\g$). Denote by
$\omega$ the Chevalley involution in $U\g$. Then the map
$\theta:x\mapsto\omega(\overline{x})$ is an isomorphism $U\n_-\to
U\n_+$, and $\mathbb S_\lambda(x\otimes
y)=\pi_\lambda(\theta(x)\otimes y)=(\omega(x)y)_0(\lambda)$ is the
Shapovalov form on $U\n_-$.

Set
\[
K_\lambda=\{y\in U\n_-\,|\,\pi_\lambda(x\otimes y)=0\ \mbox{\rm
for all}\ x\in U\n_+\},
\]
\[
\widetilde{K}_\lambda=\{x\in U\n_+\,|\,\pi_\lambda(x\otimes y)=0\
\mbox{\rm for all}\ y\in U\n_-\}.
\]
Clearly, $K_\lambda$ is the kernel of $\mathbb S_\lambda$,
$\widetilde{K}_\lambda=\omega(\overline{K_\lambda})$. Notice also
that $K(\lambda)=K_\lambda\cdot\mathbf1_\lambda$.

For any $\g$-module $L$ and subset $P\subset U\g$ define
\[
L[0]^P=\{l\in L[0]\,|\,al=\varepsilon(a)l \ \mbox{for all}\ a\in
P\}
\]
(here $\varepsilon$ stands for the standard counit in $U\g$). In
particular,
\[
F[0]^P=\{f\in F[0]\,|\,\overrightarrow{a}f=\varepsilon(a)f \
\mbox{for all}\ a\in P\}
\]

\begin{theorem}\label{thm-theta-image}
The image of $\Theta$ is $F[0]^{K_\lambda+\widetilde{K}_\lambda}$.
\end{theorem}


In order to prove Theorems \ref{thm-theta-inj} and
\ref{thm-theta-image} we need some preparations.

In the sequel $L$ stands for a $\g$-module which is a direct sum
of finite dimensional $\g$-modules.

For a $\g$-module $M$ which is a direct sum of finite-dimensional
$\h$-weight spaces we will denote by $M^*$ its restricted dual.

Let $\widetilde{M}(\lambda)$ be the ``opposite Verma module'' with
the lowest weight $\lambda\in\h^*$ and the lowest weight vector
$\widetilde{{\mathbf 1}}_\lambda$. It is clear that
$\widetilde{K}_{-\lambda}\cdot\widetilde{{\mathbf 1}}_\lambda$ is
the maximal
$\g$-submodule in $\widetilde{M}(\lambda)$.

\begin{lemma}\label{lem-n+-inv}
$\Hom_{\n_-}(M(\lambda),L)=(M(\lambda)^*\otimes L)^{\n_-}$,
$\Hom_{\n_+}(\widetilde{M}(\lambda),L)=(\widetilde{M}(\lambda)^*\otimes
L)^{\n_+}$.
\end{lemma}

\begin{proof}
For any $\varphi\in\Hom_{\n_-}(M(\lambda),L)$ the image of
$\varphi$ is equal to the finite-dimensional $\n_-$-submodule
$U\n_-\cdot\varphi({\mathbf 1}_\lambda)$. Therefore for any $x\in
U\n_-$ such that $x{\mathbf 1}_\lambda$ is a weight vector whose
weight is large enough we have $\varphi(x{\mathbf
1}_\lambda)=x\varphi({\mathbf 1}_\lambda)=0$. Thus $\varphi$
corresponds to an element in $(M(\lambda)^*\otimes L)^{\n_-}$.

The second part of the lemma can be proved similarly.
\end{proof}

Choose vectors ${\mathbf1}_{\lambda}^*\in M(\lambda)^*[-\lambda]$
and $\widetilde{{\mathbf1}}_{-\lambda}^*\in
\widetilde{M}(-\lambda)^*[\lambda]$ such that
$\langle{\mathbf1}_{\lambda}^*,{\mathbf1}_{\lambda}\rangle=
\langle\widetilde{{\mathbf1}}_{-\lambda}^*,\widetilde{{\mathbf1}}_{-\lambda}\rangle=1$.
Define maps
$\zeta:\Hom_\g(M(\lambda),\widetilde{M}(-\lambda)^*\otimes L)\to
L[0]$ and
$\widetilde{\zeta}:\Hom_\g(\widetilde{M}(-\lambda),M(\lambda)^*\otimes
L)\to L[0]$ by the formulas
$\varphi({\mathbf1}_{\lambda})=\widetilde{{\mathbf1}}_{-\lambda}^*\otimes\zeta_\varphi\,+$
lower order terms,
$\varphi(\widetilde{{\mathbf1}}_{-\lambda})=
{\mathbf1}_{\lambda}^*\otimes\widetilde{\zeta}_\varphi\,+$ higher
order terms.

Consider also the natural maps
\begin{gather*}
r:\Hom_\g(M(\lambda)\otimes
\widetilde{M}(-\lambda),L)\to\Hom_\g(M(\lambda),\widetilde{M}(-\lambda)^*\otimes
L),\\
\widetilde{r}:\Hom_\g(M(\lambda)\otimes
\widetilde{M}(-\lambda),L)\to\Hom_\g(\widetilde{M}(-\lambda),M(\lambda)^*\otimes
L).
\end{gather*}

\begin{proposition}\label{prop-zeta-r-isomorphisms}
Maps $\zeta$, $\widetilde{\zeta}$, $r$, and $\widetilde{r}$ are
vector space isomorphisms, and the diagram
\[
 \xymatrix{\Hom_\g(M(\lambda)\otimes
\widetilde{M}(-\lambda),L) \ar[d]^{\widetilde{r}} \ar[r]^{r} &
\Hom_\g(M(\lambda),\widetilde{M}(-\lambda)^*\otimes
L) \ar[d]^{\zeta} \\
\Hom_\g(\widetilde{M}(-\lambda),M(\lambda)^*\otimes L)
\ar[r]^{\widetilde{\zeta}
} & L[0]}
\]
is commutative.
\end{proposition}

\begin{proof}
First of all notice that we have the natural
identifications
\begin{gather*}
\Hom_\g(M(\lambda),\widetilde{M}(-\lambda)^*\otimes
L)=(\widetilde{M}(-\lambda)^*\otimes L)^{\n_+}[\lambda],\\
\Hom_\g(\widetilde{M}(-\lambda),M(\lambda)^*\otimes
L)=(M(\lambda)^*\otimes L)^{\n_-}[-\lambda].
\end{gather*}

Further on, we have
\begin{gather*}
\Hom_\g(M(\lambda)\otimes
\widetilde{M}(-\lambda),L)=\Hom_\g(M(\lambda),\Hom_\g(\widetilde{M}(-\lambda),L))=\\
\Hom_{\n_+}(\widetilde{M}(-\lambda),L)[\lambda]=L[0].
\end{gather*}
On the other side,
$\Hom_{\n_+}(\widetilde{M}(-\lambda),L)[\lambda]=(\widetilde{M}(-\lambda)^*\otimes
L)^{\n_+}[\lambda]$ by Lemma \ref{lem-n+-inv}. Now it is easy to
see that the map $r$ (resp.\ $\zeta$) corresponds to the
identification $\Hom_\g(M(\lambda)\otimes
\widetilde{M}(-\lambda),L)=(\widetilde{M}(-\lambda)^*\otimes
L)^{\n_+}[\lambda]$ (resp.\ $(\widetilde{M}(-\lambda)^*\otimes
L)^{\n_+}[\lambda]=L[0]$).

The second part of the proposition concerning $\widetilde{r}$ and
$\widetilde{\zeta}$ can be verified similarly.
\end{proof}

Now note that the pairing $\pi_\lambda:U\n_+\otimes
U\n_-\to\mathbb C$ naturally defines the pairing
$\widetilde{M}(-\lambda)\otimes M(\lambda)\to\mathbb C$. Denote by
$\chi_\lambda:M(\lambda)\to\widetilde{M}(-\lambda)^*$ the
corresponding morphism of $\g$-modules. The kernel of
$\chi_\lambda$ is equal to
$K(\lambda)=K_\lambda\cdot\mathbf1_\lambda$, and the image of
$\chi_\lambda$ is
$(\widetilde{K}_\lambda\cdot\widetilde{\mathbf1}_{-\lambda})^\perp$.
Therefore,
$(\widetilde{K}_\lambda\cdot\widetilde{\mathbf1}_{-\lambda})^\perp\simeq
V(\lambda)$, and $\chi_\lambda$ can be naturally represented as
$\chi_\lambda''\circ\chi_\lambda'$, where
\[
M(\lambda)\stackrel{\chi_\lambda'}{\longrightarrow}V(\lambda)
\stackrel{\chi_\lambda''}{\longrightarrow}\widetilde{M}(-\lambda)^*.
\]

The morphisms $\chi_\lambda'$ and $\chi_\lambda''$ induce the
commutative diagram of inclusions
\[
\xymatrix{\Hom_\g(V(\lambda),V(\lambda)\otimes L) \ar[d] \ar[r] &
\Hom_\g(M(\lambda),V(\lambda)\otimes L) \ar[d]\\
\Hom_\g(V(\lambda),\widetilde{M}(-\lambda)^*\otimes L) \ar[r] &
\Hom_\g(M(\lambda),\widetilde{M}(-\lambda)^*\otimes L).}
\]
It is clear that the following lemma holds:

\begin{lemma}\label{lem-image-chi-lambda-pr}
The image of $\Hom_\g(V(\lambda),V(\lambda)\otimes L)$ in
$\Hom_\g(M(\lambda),\widetilde{M}(-\lambda)^*\otimes L)$ under the
inclusion above consists of the morphisms
$\varphi:M(\lambda)\to\widetilde{M}(-\lambda)^*\otimes L$ such
that $\varphi(K_\lambda{\mathbf1}_{\lambda})=0$ and
$\varphi(M(\lambda))\subset
(\widetilde{K}_\lambda\widetilde{{\mathbf1}}_{-\lambda})^\perp\otimes
L$.\qed
\end{lemma}

\begin{proposition}\label{prop-tilde-K-lambda-0}
Let $\varphi\in\Hom_\g(M(\lambda),\widetilde{M}(-\lambda)^*\otimes
L)$. Then $\varphi(M(\lambda))\subset
(\widetilde{K}_\lambda\widetilde{{\mathbf1}}_{-\lambda})^\perp\otimes
L$ iff $\widetilde{K}_\lambda\zeta_\varphi=0$.
\end{proposition}

\begin{proof}
First notice that $\varphi(M(\lambda))\subset
(\widetilde{K}_\lambda\widetilde{{\mathbf1}}_{-\lambda})^\perp\otimes
L$ iff $\varphi(\mathbf1_\lambda)\in
(\widetilde{K}_\lambda\widetilde{{\mathbf1}}_{-\lambda})^\perp\otimes
L$. Indeed, for any $x\in U\g$ we have
$\varphi(x\mathbf1_\lambda)=\sum_{(x)}(x_{(1)}\otimes
x_{(2)})\varphi(\mathbf1_\lambda)$ and
$U\g\cdot(\widetilde{K}_\lambda\widetilde{{\mathbf1}}_{-\lambda})^\perp=
(\widetilde{K}_\lambda\widetilde{{\mathbf1}}_{-\lambda})^\perp$.

Denote by $\psi$ the element in
$\Hom_{\n_+}(\widetilde{M}(-\lambda),L)$ that corresponds to
$\varphi(\mathbf1_\lambda)\in(\widetilde{M}(-\lambda)^*\otimes
L)^{\n_+}$ (see Lemma \ref{lem-n+-inv}). Under this notation
$\varphi(\mathbf1_\lambda)\in
(\widetilde{K}_\lambda\widetilde{{\mathbf1}}_{-\lambda})^\perp\otimes
L$ iff
$\psi(\widetilde{K}_\lambda\widetilde{{\mathbf1}}_{-\lambda})=0$.
On the other hand,
$\zeta_\varphi=\psi(\widetilde{{\mathbf1}}_{-\lambda})$ and
$\psi(\widetilde{K}_\lambda\widetilde{{\mathbf1}}_{-\lambda})=
\widetilde{K}_\lambda\psi(\widetilde{{\mathbf1}}_{-\lambda})=\widetilde{K}_\lambda\zeta_\varphi$.
This completes the proof.
\end{proof}

\begin{proposition}\label{prop-K-lambda-0}
Let $\varphi\in\Hom_\g(M(\lambda),\widetilde{M}(-\lambda)^*\otimes
L)$. Then $\varphi(K_\lambda{\mathbf1}_{\lambda})=0$ iff
$K_\lambda\zeta_\varphi=0$.
\end{proposition}

\begin{proof}
Consider
$\widehat{\varphi}=r^{-1}(\varphi):\Hom_\g(M(\lambda)\otimes
\widetilde{M}(-\lambda),L)$ and
$\widetilde{\varphi}=\widetilde{r}(\widehat{\varphi})=
\widetilde{r}(r^{-1}(\varphi)):\Hom_\g(\widetilde{M}(-\lambda),M(\lambda)^*\otimes
L)$ (see Proposition \ref{prop-zeta-r-isomorphisms}). Clearly,
$\varphi(K_\lambda{\mathbf1}_{\lambda})=0$ iff
$\widehat{\varphi}(K_\lambda{\mathbf1}_{\lambda}\otimes\widetilde{M}(-\lambda))=0$
iff $\widetilde{\varphi}(\widetilde{M}(-\lambda))
\subset(K_\lambda{\mathbf1}_{\lambda})^\perp\otimes L$.

Arguing as in the proof of Proposition \ref{prop-tilde-K-lambda-0}
we see that $\widetilde{\varphi}(\widetilde{M}(-\lambda))
\subset(K_\lambda{\mathbf1}_{\lambda})^\perp\otimes L$ iff
$K_\lambda\widetilde{\zeta}_{\widetilde{\varphi}}=0$. Now it is
enough to notice that $\widetilde{\zeta}_{\widetilde{\varphi}}=
\widehat{\varphi}({\mathbf1}_{\lambda}\otimes{\widetilde{\mathbf1}}_{-\lambda})=\zeta_\varphi$.
\end{proof}

Define a map $u:\Hom_\g(V(\lambda),V(\lambda)\otimes L)\to L[0]$
via $\varphi\mapsto u_\varphi$, where
$\varphi(\overline{\mathbf1}_\lambda)=\overline{\mathbf1}_\lambda\otimes
u_\varphi\,+$ lower order terms.

\begin{proposition}
The map $u$ defines the isomorphism
$\Hom_\g(V(\lambda),V(\lambda)\otimes L)\simeq
L[0]^{K_\lambda+\widetilde{K}_\lambda}$.
\end{proposition}

\begin{proof}
Observe that $u$ can be decomposed as
\[
\Hom_\g(V(\lambda),V(\lambda)\otimes L)\longrightarrow
\Hom_\g(M(\lambda),\widetilde{M}(-\lambda)^*\otimes L)
\stackrel{\zeta}{\longrightarrow}L[0],
\]
where the first arrow is the natural inclusion considered in Lemma
\ref{lem-image-chi-lambda-pr}. Now it is enough to apply the above
mentioned lemma and Propositions \ref{prop-tilde-K-lambda-0} and
\ref{prop-K-lambda-0}.
\end{proof}

Applying the last proposition to the case $L=F$ we get Theorems
\ref{thm-theta-inj} and \ref{thm-theta-image}.

\medskip

Now we describe $\Theta^{-1}:
F[0]^{K_\lambda+\widetilde{K}_\lambda}\to\Hom_\g(V(\lambda),V(\lambda)\otimes
F)$ explicitly. We are going to obtain a formula for
$\Theta^{-1}$ in terms of the Shapovalov form.
By means of the standard identification $M(\lambda)\simeq U\n_-$
we can regard $\mathbb S_\lambda$ as a bilinear form on
$M(\lambda)$. Denote by $\overline{\mathbb S}_\lambda$ the
corresponding bilinear form on $V(\lambda)$. Set
\[
Q_+=\left(\sum_{\alpha\in\Pi}\mathbb
Z_+\alpha\right)\setminus\{0\}.
\]
For any
$\beta\in Q_+$ denote by
$\overline{\mathbb S}_\lambda^\beta$ the restriction of
$\overline{\mathbb S}_\lambda$ to $V(\lambda)[\lambda-\beta]$.
Let
$x_\beta^i\cdot\overline{\mathbf1}_\lambda$ be an arbitrary basis
in $V(\lambda)[\lambda-\beta]$, where $x_\beta^i\in
U\n_-[-\beta]$.

Take $f\in F[0]^{K_\lambda+\widetilde{K}_\lambda}$ and set
$\varphi=\Theta^{-1}(f)$, i.e.,
\[
\varphi(\overline{\mathbf1}_\lambda)=\overline{\mathbf1}_\lambda\otimes
f+\sum_{\beta\in
Q_+}\sum_ix_\beta^i\cdot\overline{\mathbf1}_\lambda\otimes
f^{\beta,i}.
\]

\begin{proposition}\label{prop-Theta-inv}
$f^{\beta,i}=\sum_j\left(\overline{\mathbb
S}_\lambda^\beta\right)^{-1}_{ij}\overrightarrow{\omega\left(\overline{x_\beta^j}\right)}f$.
\end{proposition}

\begin{proof}
Set $\xi=\varphi(\overline{\mathbf1}_\lambda)$. Clearly, $\xi$ is
a singular element in $V(\lambda)\otimes F$. In particular,
$(e \otimes 1 + 1 \otimes e)\xi = 0$, i.e., $(e \otimes
1)\xi=(1 \otimes \overline{e})\xi$ for any $e\in\n_+$. By
induction we get $(x \otimes 1)\xi=(1 \otimes
\overline{x})\xi$ for any
$x\in U\n_+$. Therefore
\[
\left(\overline{\mathbb
S}_\lambda\otimes\id\right)
\left(\overline{\mathbf1}_\lambda\otimes
\left(\omega\left(x_\beta^j\right)\otimes1\right)\xi\right)
=\left(\overline{\mathbb S}_\lambda\otimes\id\right)
\left(\overline{\mathbf1}_\lambda\otimes\left(1\otimes
\omega\left(\overline{x_\beta^j}\right)\right)\xi\right).
\]
Calculating both sides of this equation we get
\[
\sum_i{\mathbb S}_\lambda(x_\beta^j\otimes x_\beta^i)f^{\beta,i}=
\overrightarrow{\omega\left(\overline{x_\beta^j}\right)}f,
\]
and the proposition follows.
\end{proof}

\medskip

Let us define an associative product $\star_\lambda$ on
$F[0]^{K_\lambda+\widetilde{K}_\lambda}$ by means of $\Theta$. We
are going to obtain an explicit formula for $\star_\lambda$ in
terms of the Shapovalov form.

\begin{theorem}
For any $f_1,f_2\in F[0]^{K_\lambda+\widetilde{K}_\lambda}$ we
have
\begin{equation}\label{sym_space_star_prod}
f_1\star_\lambda
f_2=\mu\left(\overrightarrow{J_{\red}(\lambda)}(f_1\otimes
f_2)\right),
\end{equation}
where
\begin{equation}\label{J_red}
J_{\red}(\lambda)=1\otimes1+\sum_{\beta\in
Q_+}\sum_{i,j}\left(\overline{\mathbb
S}_\lambda^\beta\right)^{-1}_{ij}
x_\beta^i\otimes\omega\left(\overline{x_\beta^j}\right).
\end{equation}
\end{theorem}

\begin{proof}
We have $f_1\star_\lambda f_2=\Theta(\varphi_1\ast\varphi_2)$,
where $\varphi_1=\Theta^{-1}(f_1)$, $\varphi_2=\Theta^{-1}(f_2)$,
and $\ast$ is the product on $\Hom_\g(V(\lambda),V(\lambda)\otimes
F)$ given by \eqref{eqn-basic-star-prod}. Now observe that
\begin{gather*}
(\varphi_1\ast\varphi_2)(\overline{\mathbf1}_\lambda)=
(\id\otimes\mu)(\varphi_1\otimes\id)(\varphi_2(\overline{\mathbf1}_\lambda))=\\
(\id\otimes\mu)(\varphi_1\otimes\id)\left(\overline{\mathbf1}_\lambda\otimes
f_2+\sum_{\beta\in
Q_+}\sum_ix_\beta^i\cdot\overline{\mathbf1}_\lambda\otimes
f_2^{\beta,i}\right)=\\
(\id\otimes\mu)\left(\varphi_1(\overline{\mathbf1}_\lambda)\otimes
f_2+\sum_{\beta\in
Q_+}\sum_i(\Delta(x_\beta^i)\varphi_1(\overline{\mathbf1}_\lambda))\otimes
f_2^{\beta,i}\right)=\\
\overline{\mathbf1}_\lambda\otimes\left(f_1f_2+\sum_{\beta\in Q_+}
\sum_i\left(\overrightarrow{x_\beta^i}f_1\right)f_2^{\beta,i}\right)+\mbox{lower
order terms}.
\end{gather*}
Therefore
\[
f_1\star_\lambda f_2=f_1f_2+\sum_{\beta\in
Q_+}\sum_i\left(\overrightarrow{x_\beta^i}f_1\right)f_2^{\beta,i}.
\]
To finish the proof it is enough now to apply Proposition
\ref{prop-Theta-inv} to $f_2$.
\end{proof}



\subsection{Application to symmetric spaces}\label{Subsect-main-symm}

Now let us apply the construction above to some specific values of
$\lambda\in\h^*$.

Let $\Delta\subset\Pi$. Assume that $\lambda\in\h^*$ is such that
$\langle\lambda,\alpha^\vee\rangle=n_\alpha\in\mathbb Z_+$ for any
$\alpha\in\Delta$, and
$\langle\lambda+\rho,\beta^\vee\rangle\not\in\mathbb N$ for $\beta\in\B
R_+\setminus\spanv\Delta$.

\begin{proposition}
Let $L$ be a $\g$-module which is a direct sum of
finite-dimen\-si\-o\-nal $\g$-modules, $l\in L[0]$. Then $l\in
L[0]^{K_\lambda}$ iff $l\in L[0]^{\widetilde{K}_\lambda}$.
\end{proposition}

\begin{proof}
Recall that in our case $K_\lambda$ (resp.\
${\widetilde{K}_\lambda}$) is generated by
$Y_\alpha^{n_\alpha+1}$ (resp.\
$X_\alpha^{n_\alpha+1}$) for all
$\alpha\in\Delta$.

Now take any $\alpha\in\Delta$ and regard $L$ as an
$\sl(2)_\alpha$-module, where
$\sl(2)_\alpha\subset\g$ is a subalgebra generated by $X_\alpha$
and $Y_\alpha$. By standard structure theory of finite-dimensional
$\sl(2)$-modules we see that $Y_\alpha^{n_\alpha+1}l=0$ iff
$X_\alpha^{n_\alpha+1}l=0$, which completes the proof.
\end{proof}

\begin{corollary}
$F[0]^{K_\lambda+\widetilde{K}_\lambda}=F[0]^{K_\lambda}$.\qed
\end{corollary}

Therefore in this case we get the associative algebra
$(F[0]^{K_\lambda},\star_\lambda)$.

Now assume additionally that $n_\alpha=0$ for all
$\alpha\in\Delta$. In this case $F[0]^{K_\lambda}$ is closed under
the original product in $F$. Moreover,
$F[0]^{K_\lambda}\simeq\Fun(G/U)$, where $U\subset G$ is the
reductive Levi subgroup that corresponds to $\Delta$. Notice that
in this case formula \eqref{sym_space_star_prod} defines an
equivariant quantization of the Kirillov-Kostant-Souriau bracket
on the coadjoint orbit through $\lambda$. A formula of this type
appears also in \cite{AL}.

\begin{remark}\label{rem-Mazorchuk}
Assume again that $\lambda\in\h^*$ is as decribed in the beginning
of this subsection ($n_\alpha$ should not necessary be $0$). It is
known that in this case $M(\lambda)$ is projective (see
\cite{Ja}).
In particular, the natural map $\Hom_\g(M(\lambda),
M(\lambda)\otimes F)\to\Hom_\g(M(\lambda), V(\lambda)\otimes F)$
is surjective. This map factors as $\Hom_\g(M(\lambda),
M(\lambda)\otimes F)\to\Hom_\g(V(\lambda), V(\lambda)\otimes
F)\to\Hom_\g(M(\lambda), V(\lambda)\otimes F)$, where the first
map defined by Proposition \ref{prop-Hopf-submod}. Hence the map
$\Hom_\g(M(\lambda), M(\lambda)\otimes F)\to\Hom_\g(V(\lambda),
V(\lambda)\otimes F)$ is also surjective. By Proposition
\ref{prop-Hom-End-equiv} the canonical map $(\End
M(\lambda))^{\rrr}_{\fin} \to (\End V(\lambda))^{\rrr}_{\fin}$ is
also surjective. Since the natural map $U\g\to(\End
M(\lambda))^{\rrr}_{\fin}$ is surjective for any $\lambda$, we
recover the fact that the map $U\g\to\End
V(\lambda))^{\rrr}_{\fin}$ is also surjective in this case (see
\cite{Jo}).
\end{remark}

\section{Limiting properties of the fusion element}\label{Sect-limit}

In this section we will freely use the notation of the previous
one. For any generic $\lambda\in\h^*$ (i.e.,
$\langle\lambda_0+\rho,\beta^\vee\rangle\not\in\mathbb N$ for all
$\beta\in\B R_+$) we denote by $J(\lambda)$ the fusion element
related to the Verma module $M(\lambda)$ (see, e.g.,
\cite{ES_DYBE}). Notice that in this case $V(\lambda)=M(\lambda)$
and $J_{\red}(\lambda)=J(\lambda)$.

\subsection{One distinguished root case}\label{Subsect-limit-one-root}

Fix $\alpha\in\B R_+$. Take $\lambda_0\in\h^*$ such that
$\langle\lambda_0+\rho,\alpha^\vee\rangle=n\in\mathbb N$,
$\langle\lambda_0+\rho,\beta^\vee\rangle\not\in\mathbb N$ for all
$\beta\in\B R_+\setminus\{\alpha\}$.

\begin{theorem}\label{thm-reg-one-root-1}
Let $N$ be an arbitrary $\n_+$-module. Consider the family of
operators $J(\lambda)_N:F[0]^{K_{\lambda_0}}\otimes N\to F\otimes
N$ naturally defined by $J(\lambda)$. Then this family is regular
at $\lambda=\lambda_0$.
\end{theorem}

\begin{proof}
Fix an arbitrary line $l\subset\h^*$ through $\lambda_0$,
$l=\{\lambda_0+t\nu\,|\,t\in\mathbb C\}$, transversal to the
hyperplane $\langle\lambda+\rho,\alpha^\vee\rangle=n$.

Identify $M(\lambda)$ with $U\n_-$ in the standard way. Recall
that we have a basis $x_\beta^i\in U\n_-[-\beta]$ for $\beta\in
Q_+$. Let $L\left(\mathbb S_\lambda^\beta\right)\in\End
U\n_-[-\beta]$ be given by the matrix $\left(\mathbb
S_\lambda^\beta\right)_{ij}$ in the basis
$x_\beta^i$. Notice that $\Ker L\left(\mathbb
S_{\lambda_0}^\beta\right)=\Ker\mathbb
S_{\lambda_0}^\beta=K_{\lambda_0}[-\beta]:=K_{\lambda_0}\cap
U\n_-[-\beta]$. For any $\lambda\in l$ sufficiently close to
$\lambda_0$, $\lambda\neq\lambda_0$ we have $M(\lambda)$ is
irreducible, and $L\left(\mathbb S_\lambda^\beta\right)$ is
invertible for any $\beta\in Q_+$. In this notation we have
\[
J(\lambda)=1\otimes1+\sum_{\beta\in Q_+}\sum_j L\left(\mathbb
S_\lambda^\beta\right)^{-1}x_\beta^j\otimes\omega\left(\overline{x_\beta^j}\right).
\]

Take $\lambda=\lambda_0+t\nu\in l$. Fix any $\beta\in Q_+$ and set
$V=U\n_-[-\beta]$, $A_t=L\left(\mathbb S_\lambda^\beta\right)$,
$V_0=\Ker A_0=K_{\lambda_0}[-\beta]\subset V$. Write
$A_t=A_0+tB_t$, where $B_t$ is regular at $t=0$. It is known (see,
e.g., \cite{ESt_1}) that we have $A_t^{-1}=\frac{1}{t}C+D_t$,
where
$D_t$ is regular at $t=0$.

\begin{lemma}\label{lem-limit-lin-alg-1}
$\Image C\subset V_0$.
\end{lemma}

\begin{proof}
We have $A_tA_t^{-1}=\id$ for any $t\neq0$, i.e.,
$\frac{1}{t}A_0C+A_0D_t+B_tC+tB_tD_t=\id$. Since the left hand
side should be regular at $t=0$, we have $A_0C=0$, which proves
the lemma.
\end{proof}

For $t\neq 0$ set
$J_t=\sum_jA_t^{-1}x_j\otimes\omega(\overline{x_j})$ (from now on
we are omitting the index $\beta$ for the sake of brevity). By
Lemma \ref{lem-limit-lin-alg-1} we have
$Cx_j\in V_0=K_{\lambda_0}[-\beta]$. Hence for
$f\in F[0]^{K_{\lambda_0}}$ we have $\overrightarrow{Cx_j}f=0$. Therefore
$\overrightarrow{A_t^{-1}x_j}f=\frac{1}{t}\overrightarrow{Cx_j}f+\overrightarrow{D_tx_j}f=
\overrightarrow{D_tx_j}f$. This proves the regularity of
$(J_t)_N(f\otimes \cdot)$ at $t=0$, i.e., the regularity
of $J(\lambda)_N(f\otimes \cdot)$ at
$\lambda=\lambda_0$.
\end{proof}

Similarly to Theorem \ref{thm-reg-one-root-1} one can prove the
following

\begin{theorem}\label{thm-reg-one-root-2}
Let $M$ be an arbitrary $\n_-$-module. Consider the family of
operators $J(\lambda)^M:M\otimes
F[0]^{\widetilde{K}_{\lambda_0}}\to M\otimes F$ naturally defined
by $J(\lambda)$. Then this family is regular at
$\lambda=\lambda_0$.\qed
\end{theorem}

\begin{theorem}\label{thm-limit-one-root}
Let $f\in F[0]^{K_{\lambda_0}}$, $g\in
F[0]^{\widetilde{K}_{\lambda_0}}$. Then
$\overrightarrow{J(\lambda)}(f\otimes g)\rightarrow
\overrightarrow{J_{\red}(\lambda_0)}(f\otimes g)$ as $\lambda\rightarrow\lambda_0$.
\end{theorem}

\begin{proof}
We will use the notation defined in the proof of Theorem
\ref{thm-reg-one-root-1}.

\begin{lemma}\label{lem-limit-lin-alg-2}
$D_0A_0=\id$ on $V/V_0$.
\end{lemma}

\begin{proof}
Arguing in the same manner as in the proof of Lemma
\ref{lem-limit-lin-alg-1} but starting from
$A_t^{-1}A_t=\id$ and setting $t=0$, we get $D_0A_0+CB_0=\id$. Now
notice that for any $v\in V$ we have $CB_0v\in V_0$, which proves
the lemma.
\end{proof}

Since the Shapovalov form is symmetric, we may choose
$V_1\subset V$ such that $V=V_0\oplus V_1$, $A_0(V_1)=V_1$, and $A_0$ is non-degenerate on $V_1$.
Assume that the basis $x_j$ is compatible with this decomposition.
We see that
$\overrightarrow{A_t^{-1}x_j}f\rightarrow
\overrightarrow{D_0x_j}f$ as $t\to0$.

For any $x_j\in V_0$ we have
$\omega(\overline{x_j})\in\widetilde{K}_{\lambda_0}\cap
U\n_+[\beta]$. This implies that
$\overrightarrow{D_0x_j}f\otimes\overrightarrow{\omega(\overline{x_j})}g=0$
by our assumptions on $g$.

For any $x_j\in V_1$ we see, by Lemma
\ref{lem-limit-lin-alg-2}, that
$\overrightarrow{D_0x_j}f=\overrightarrow{A_0^{-1}x_j}f$. Thus
\[
\overrightarrow{J_t}(f\otimes g)\rightarrow
\sum_{j: x_j\in
V_1}\overrightarrow{A_0^{-1}x_j}f\otimes\overrightarrow{\omega(\overline{x_j})}g.
\]
Clearly, this means, by definition of $J_{\red}(\lambda_0)$, that
$\overrightarrow{J(\lambda)}(f\otimes g)\rightarrow
\overrightarrow{J_{\red}(\lambda_0)}(f\otimes g)$ as
$\lambda\rightarrow\lambda_0$.
\end{proof}

\begin{corollary}\label{cor-limit-one-root}
Let $f_1,f_2\in F[0]^{K_{\lambda_0}+\widetilde{K}_{\lambda_0}}$.
Then $f_1\star_\lambda f_2\rightarrow f_1\star_{\lambda_0}f_2$ as
$\lambda\rightarrow\lambda_0$.\qed
\end{corollary}

\begin{example}
Let $\g=\sl(2)$. In \cite{EKKMASVT} we considered the star-product
on polynomial functions on coadjoint orbits $\mathcal O_\lambda$
of $\g$ defined by the natural action of the fusion element
$J(\lambda)$. In particular, we obtained the formula
\[
f_a\star_\lambda
f_b=\left(1-\frac{1}{\lambda}\right)f_af_b+\frac{1}{2}f_{[a,b]}
+\frac{\lambda}{2}\langle a,b\rangle,
\]
where $f_x$ is the restriction onto $\mathcal O_\lambda$ of the
linear function on $\g^*$ defined by $x\in\g$, and $\langle
a,b\rangle=\Tr(ab)$. Despite $J(\lambda)$ has a singularity at
$\lambda=1$ we see that $f_a\star_1 f_b$ is well defined, and the
set $\{f_x\,|\,x\in\g\}$ generates an algebra under $\star_1$
isomorphic to $\End V(1)\simeq \Mat(2,\mathbb C)$.

Similarly, one can also show that for any $\lambda\in\mathbb Z_+$
the set $\{f_x\,|\,x\in\g\}$ generates an algebra under
$\star_\lambda$ isomorphic to $\End V(\lambda)\simeq \Mat(\lambda+1,\mathbb
C)$. Corollary \ref{cor-limit-one-root} explains these phenomena.
\end{example}

\subsection{Regularity properties}\label{Subsect-regularity-general}

Let $\lambda_0\in\h^*$. We will say that $\lambda_0$ has the {\it
good regularity property} if for any $\n_-$-module $M$ the family
of operators $J(\lambda)^M:M\otimes
F[0]^{\widetilde{K}_{\lambda_0}}\to M\otimes F$ naturally defined
by $J(\lambda)$ is regular at $\lambda=\lambda_0$. Clearly, if
$\lambda_0$ is generic (i.e., $V(\lambda_0)=M(\lambda_0)$ is
irreducible), then $\lambda_0$ has the good regularity property.
We have seen that $\lambda_0$ as in Subsection
\ref{Subsect-limit-one-root} also has the good regularity
property.

\begin{theorem}\label{thm-limit-general}
Assume that $\lambda_0\in\h^*$ has the good regularity property.
Then for any $f\in F[0]^{K_{\lambda_0}}$, $g\in
F[0]^{\widetilde{K}_{\lambda_0}}$ we have
$\overrightarrow{J(\lambda)}(f\otimes g)\rightarrow
\overrightarrow{J_{\red}(\lambda_0)}(f\otimes g)$ as
$\lambda\rightarrow\lambda_0$.
\end{theorem}

\begin{proof}
For any $\lambda\in\h^*$ we may naturally identify $M(\lambda)$
with $U\n_-$ as $\n_-$-modules. Therefore we know by definition of
a good regular property that
$J(\lambda)^{M(\lambda)}(\mathbf1_\lambda\otimes g)$ is regular at
$\lambda=\lambda_0$. Thus
$J(\lambda)^{M(\lambda)}(\mathbf1_\lambda\otimes g)\rightarrow
Z\in M(\lambda_0)\otimes F$ as $\lambda\rightarrow\lambda_0$. In
an arbitrary basis $x_\beta^i\in U\n_-[-\beta]$ we have
\[
J(\lambda)^{M(\lambda)}(\mathbf1_\lambda\otimes
g)=\mathbf1_\lambda\otimes g\,+\sum_{\beta\in
Q_+}\sum_{i,j}\left(\mathbb S_\lambda^\beta\right)^{-1}_{ij}
x_\beta^i\mathbf1_\lambda\otimes\overrightarrow{\omega\left(\overline{x_\beta^j}\right)}g,
\]
and
\[
Z=\mathbf1_{\lambda_0}\otimes g+\sum_{\beta\in
Q_+}\sum_{i,j}a^\beta_{ij}
\left(x_\beta^i\mathbf1_{\lambda_0}\right)\otimes
\overrightarrow{\omega\left(\overline{x_\beta^j}\right)}g
\]
for some coefficients $a^\beta_{ij}\in\mathbb C$.

Now choose a basis $x_\beta^i\in U\n_-[-\beta]$ in the following
way: first take a basis in
$K_{\lambda_0}[-\beta]=K_{\lambda_0}\cap U\n_-[-\beta]$ and then
extend it arbitrarily to a basis in the whole $U\n_-[-\beta]$. In
this basis the projection $\overline{Z}\in V(\lambda_0)\otimes F$
of the element $Z$ is given by
\begin{equation}\label{eqn-Z-bar}
\overline{Z}=\overline{\mathbf1}_{\lambda_0}\otimes g+\sum_{\beta\in
Q_+}\sum_{\ x_\beta^i,x_\beta^j\not\in
K_{\lambda_0}[-\beta]}a^\beta_{ij}
\left(x_\beta^i\overline{\mathbf1}_{\lambda_0}\right)\otimes
\overrightarrow{\omega\left(\overline{x_\beta^j}\right)}g.
\end{equation}

Now notice that $Z$, being the limit of singular vectors of weight
$\lambda$ in
$M(\lambda)\otimes F$, defines the intertwining operator
$\varphi_Z\in\Hom_\g(M(\lambda_0), M(\lambda_0)\otimes F)$,
$\varphi_Z({\mathbf1}_{\lambda_0})=Z$. Under the natural map
$\Hom_\g(M(\lambda_0), M(\lambda_0)\otimes F)\to\Hom_\g(V(\lambda_0), V(\lambda_0)\otimes
F)$ we have $\varphi_Z\mapsto\varphi_{\overline{Z}}$, where
$\varphi_{\overline{Z}}(\overline{\mathbf1}_{\lambda_0})=\overline{Z}$.
Therefore
$\overline{Z}=J_{\red}(\lambda_0)^{M(\lambda_0)}(\overline{\mathbf1}_{\lambda_0}\otimes
g)$ by Proposition \ref{prop-Theta-inv} and the definition of
$J_{\red}(\lambda_0)$. Comparing this with \eqref{eqn-Z-bar} we
conclude that for all $i, j$ such that $x_\beta^i,x_\beta^j\not\in
K_{\lambda_0}[-\beta]$ we have
$a^\beta_{ij}=\left(\overline{\mathbb
S}_{\lambda_0}^\beta\right)^{-1}_{ij}$.

Finally,
\begin{gather*}
\overrightarrow{J(\lambda)}(f\otimes g)\rightarrow
fg+\sum_{\beta\in Q_+}\sum_{i,j}a^\beta_{ij}
\overrightarrow{x_\beta^i}f\otimes
\overrightarrow{\omega\left(\overline{x_\beta^j}\right)}g=\\
fg+\sum_{\beta\in Q_+}\sum_{\ x_\beta^i,x_\beta^j\not\in
K_{\lambda_0}[-\beta]}a^\beta_{ij}
\overrightarrow{x_\beta^i}f\otimes
\overrightarrow{\omega\left(\overline{x_\beta^j}\right)}g=\\
fg+\sum_{\beta\in Q_+}\sum_{\ x_\beta^i,x_\beta^j\not\in
K_{\lambda_0}[-\beta]}\left(\overline{\mathbb
S}_{\lambda_0}^\beta\right)^{-1}_{ij}
\overrightarrow{x_\beta^i}f\otimes
\overrightarrow{\omega\left(\overline{x_\beta^j}\right)}g=\\
\overrightarrow{J_{\red}(\lambda_0)}(f\otimes g)
\end{gather*}
as $\lambda\to\lambda_0$.
\end{proof}

\begin{remark}
Theorem \ref{thm-limit-general} provides another proof of Theorem
\ref{thm-limit-one-root}.
\end{remark}

\begin{corollary}\label{cor-limit-general}
Assume that $\lambda_0\in\h^*$ has the good regularity property.
Let $f_1,f_2\in F[0]^{K_{\lambda_0}+\widetilde{K}_{\lambda_0}}$.
Then $f_1\star_\lambda f_2\rightarrow f_1\star_{\lambda_0}f_2$ as
$\lambda\rightarrow\lambda_0$.\qed
\end{corollary}

\begin{proposition}
Assume that $\lambda_0\in\h^*$ has the good regularity property.
Then $F[0]^{K_{\lambda_0}}=F[0]^{\widetilde{K}_{\lambda_0}}=
F[0]^{K_{\lambda_0}+\widetilde{K}_{\lambda_0}}$.
\end{proposition}

\begin{proof}
Let $u\in F[0]^{\widetilde{K}_{\lambda_0}}$. If $\lambda\in\h^*$
is generic, then the element
$J(\lambda)^{M(\lambda)}(\mathbf1_\lambda\otimes u)$ is a singular
vector of weight $\lambda$ in $M(\lambda)\otimes F$. Therefore
$Z:=\lim_{\lambda\to\lambda_0}J(\lambda)^{M(\lambda)}(\mathbf1_\lambda\otimes
u)$ is a singular vector of weight $\lambda_0$ in
$M(\lambda_0)\otimes F$, and hence we have
$\varphi_Z\in\Hom_\g(M(\lambda_0), M(\lambda_0)\otimes F)$,
$\varphi_Z({\mathbf1}_{\lambda_0})=Z$.

Under the natural map
$\Hom_\g(M(\lambda_0), M(\lambda_0)\otimes F)\to\Hom_\g(V(\lambda_0), V(\lambda_0)\otimes
F)$ we have $\varphi_Z\mapsto\varphi_{\overline{Z}}$, where
$\varphi_{\overline{Z}}(\overline{\mathbf1}_{\lambda_0})=\overline{Z}=$ the projection
of $Z$ onto $V(\lambda_0)\otimes F$. Now notice that
$u=\Theta(\varphi_{\overline{Z}})\in
F[0]^{K_{\lambda_0}+\widetilde{K}_{\lambda_0}}$, which proves the
proposition.
\end{proof}

\begin{proposition}\label{prop-surj-1-general}
Assume that $\lambda_0\in\h^*$ has the good regularity property.
Then the natural map $\Hom_\g(M(\lambda_0), M(\lambda_0)\otimes
F)\to\Hom_\g(V(\lambda_0), V(\lambda_0)\otimes F)$ is surjective.
\end{proposition}

\begin{proof}
Recall that we have the isomorphism
\[
\Theta:\Hom_\g(V(\lambda_0),
V(\lambda_0)\otimes F)\to
F[0]^{K_{\lambda_0}+\widetilde{K}_{\lambda_0}}=F[0]^{K_{\lambda_0}}.
\]
Now take $u\in F[0]^{\widetilde{K}_{\lambda_0}}$. Consider
$Z=\lim_{\lambda\to\lambda_0}J(\lambda)^{M(\lambda)}(\mathbf1_\lambda\otimes
u)\in M(\lambda_0)\otimes F$. Since $Z$ a singular vector of
weight $\lambda_0$, we have $\varphi_Z\in\Hom_\g(M(\lambda_0),
M(\lambda_0)\otimes F)$, $\varphi_Z({\mathbf1}_{\lambda_0})=Z$.
Clearly, under the mapping $\Hom_\g(M(\lambda_0),
M(\lambda_0)\otimes F)\to\Hom_\g(V(\lambda_0), V(\lambda_0)\otimes
F)$ the image of $\varphi_Z$ equals to $\Theta^{-1}(u)$, which
proves the proposition.
\end{proof}

\begin{proposition}
Assume that $\lambda_0\in\h^*$ has the good regularity property.
Then the action map $U\g\to(\End V(\lambda_0))^{\rrr}_{\fin}$ is
surjective.
\end{proposition}

\begin{proof}
Recall that by Proposition \ref{prop-Hom-End-isom} we have the
isomorphisms
\begin{gather*}
\Hom_\g (M(\lambda_0),
M(\lambda_0)\otimes F)\simeq(\End M(\lambda_0))^{\rrr}_{\fin},\\
\Hom_\g (V(\lambda_0), V(\lambda_0)\otimes F)\simeq(\End
V(\lambda_0))^{\rrr}_{\fin}.
\end{gather*}
It is well known that the action map $U\g\to(\End
M(\lambda_0))^{\rrr}_{\fin}$ is surjective for any
$\lambda_0\in\h^*$ (see \cite{Jo}).
Since by Proposition \ref{prop-surj-1-general} the
map $(\End M(\lambda_0))^{\rrr}_{\fin}\to(\End
V(\lambda_0))^{\rrr}_{\fin}$ is surjective, the map $U\g\to(\End
V(\lambda_0))^{\rrr}_{\fin}$ is also surjective.
\end{proof}

\subsection{Symmetric space case}\label{Subsect-limit-symm}

Let $\Delta\subset\Pi$. Assume that $\lambda_0\in\h^*$ is such
that
$\langle\lambda_0,\alpha^\vee\rangle=0$ for any
$\alpha\in\Delta$, and
$\langle\lambda_0+\rho,\beta^\vee\rangle\not\in\mathbb N$ for $\beta\in\B
R_+\setminus\spanv\Delta$.

\begin{theorem}\label{thm-reg-symmetric-1}
Let $N$ be an arbitrary $\n_+$-module. Consider the family of
operators $J(\lambda)_N:F[0]^{K_{\lambda_0}}\otimes N\to F\otimes
N$ naturally defined by $J(\lambda)$. Then this family is regular
at $\lambda=\lambda_0$.
\end{theorem}

\begin{proof}
It is known (see, e.g., \cite{ESt_1}) that the only singularities
of $J(\lambda)$ near $\lambda_0$ are simple poles on the
hyperplanes
$\langle\lambda,\alpha^\vee\rangle=0$ for $\alpha\in\B
R_+\cap\spanv\Delta$. Therefore it is enough to show that for any
$f\in F[0]^{K_{\lambda_0}}$ the operator $J(\lambda)_N(f\otimes \cdot)$
has no singularity at any
such hyperplane.

Let $\Delta=\{\alpha_1,\ldots,\alpha_l\}$. For each $i=1,\ldots,l$
take an arbitrary $\lambda_i\in\h^*$ such that
$\langle\lambda_i,\alpha_i^\vee\rangle=0$, and
$\langle\lambda_i+\rho,\beta^\vee\rangle\not\in\mathbb N$ for $\beta\in\B
R_+\setminus\{\alpha_i\}$. It is well known that
$K_{\lambda_0}=K_{\lambda_1}+\ldots+K_{\lambda_l}$. In particular,
$K_{\lambda_i}\subset K_{\lambda_0}$. Also,
$F[0]^{K_{\lambda_i}}\supset F[0]^{K_{\lambda_0}}$.
Therefore we may apply
Theorem \ref{thm-reg-one-root-1} and conclude that
$J(\lambda)_N(f\otimes \cdot)$ is regular at $\lambda=\lambda_i$
for each $i$.

Now consider a hyperplane $\langle\lambda,\alpha^\vee\rangle=0$
for $\alpha\in\B R_+\cap\spanv\Delta$ which may be composite. Take
an arbitrary $\lambda'\in\h^*$ such that
$\langle\lambda',\alpha^\vee\rangle=0$, and
$\langle\lambda'+\rho,\beta^\vee\rangle\not\in\mathbb N$ for $\beta\in\B
R_+\setminus\{\alpha\}$. It follows from the results of \cite{FFM}
that $K_{\lambda'}\subset K_{\lambda_1}+\ldots+K_{\lambda_l}$,
i.e., $K_{\lambda'}\subset K_{\lambda_0}$. Arguing as above we see
that $J(\lambda)_N(f\otimes \cdot)$ is regular at
$\lambda=\lambda'$, which completes the proof.
\end{proof}

By similar considerations one can prove the following

\begin{theorem}\label{thm-reg-symmetric-2}
Let $M$ be an arbitrary $\n_-$-module. Consider the family of
operators $J(\lambda)^M:M\otimes
F[0]^{\widetilde{K}_{\lambda_0}}\to M\otimes F$ naturally defined
by $J(\lambda)$. Then this family is regular at
$\lambda=\lambda_0$.\qed
\end{theorem}

Hence we conclude that any $\lambda_0$ as described at the
beginning of this subsection has the good regularity property. In
particular, all results of Subsection
\ref{Subsect-regularity-general} are applicable to this situation.

\begin{remark}
Recall
that for $\lambda_0$ as described above the formula
$f_1\star_{\lambda_0}
f_2=\mu\left(\overrightarrow{J_{\red}(\lambda_0)}(f_1\otimes
f_2)\right)$ gives an equivariant quantization of the
Kirillov-Kostant-Souriau bracket on the coadjoint orbit through
$\lambda_0$. Applying results of Subsection
\ref{Subsect-regularity-general} we conclude that
\[
f_1\star_{\lambda_0}f_2 =
\lim_{\lambda\rightarrow\lambda_0}f_1\star_{\lambda}f_2 =
\lim_{\lambda\rightarrow\lambda_0}\mu\left(\overrightarrow{J(\lambda)}(f_1\otimes
f_2)\right).
\]
\end{remark}

\subsection{Concluding remarks}\label{Subsect-limit-open}


Let
$\Delta\subset\Pi$. It would be interesting to investigate whether our good regularity
property still holds for any $\lambda_0\in\h^*$ such that
$\langle\lambda_0,\alpha^\vee\rangle=n_\alpha\in\mathbb Z_+$ for any
$\alpha\in\Delta$, and
$\langle\lambda_0+\rho,\beta^\vee\rangle\not\in\mathbb N$ for $\beta\in\B
R_+\setminus\spanv\Delta$. This would imply that the action map
$U\g\to(\End V(\lambda_0))^{\rrr}_{\fin}$ is surjective. The
latter fact is known
for $\lambda_0$ of consideration (cf.\ also Remark
\ref{rem-Mazorchuk}). Proving the good regularity property will
provide a new explanation of this result.

\end{sloppy}

\small

\noindent E.K.: Department of Mathematics, Kharkov National University,\\
4 Svobody Sq., Kharkov \,61077, Ukraine;\\
Department of Mathematics, University of Notre Dame,\\
Notre Dame, IN \,46556, USA\\
e-mail: {\small \tt eugene.a.karolinsky@univer.kharkov.ua;
ykarolin@nd.edu}

\medskip

\noindent A.S.: Department of Mathematics, University of G\"oteborg,\\
SE-412 96 G\"oteborg, Sweden\\
e-mail: {\small \tt astolin@math.chalmers.se}

\medskip

\noindent V.T.: St.\,Petersburg Branch of Steklov Mathematical Institute,\\
Fontanka 27, St.\,Petersburg \,191023, Russia;\\
Department of Mathematical Sciences, IUPUI,\\
Indianapolis, IN 46202, USA \\
e-mail: {\small \tt vt@pdmi.ras.ru; vt@math.iupui.edu}

\end{document}